                          \def\version{7 May, 2008}                               %

\documentclass[11pt, reqno]{amsart}
\usepackage[centertags]{amsmath}
\usepackage{amsthm, a4, latexsym, amssymb}
\usepackage[latin1]{inputenc}
\usepackage{color}
\def\@rmrk#1#2{\refstepcounter
    {#1}\@ifnextchar[{\@yrmrk{#1}{#2}}{\@xrmrk{#1}{#2}}}

\makeatletter\@addtoreset{equation}{section}\makeatother

\newfont{\bfit}{cmbxti10 scaled 2000}
\newfont{\biggi}{cmr12 scaled 2000}
\newtheorem{step}{STEP}

\newcommand{\bes}{\begin{step}}
\newcommand{\es}{\end{step}}

 \newcommand{\eps}{\varepsilon}

 \newcommand{\R}{\mathbb{R}}
 \newcommand{\Z}{\mathbb{Z}}
 \newcommand{\N}{\mathbb{N}}

 \newcommand{\prob}{\mathbb{P}}

 \renewcommand{\P}{\mathbb{P}}
 \newcommand{\E}{\mathbb{E}}

 \newcommand{\heap}[2]{\genfrac{}{}{0pt}{}{#1}{#2}}

\def\1{{\mathchoice {1\mskip-4mu\mathrm l}      
{1\mskip-4mu\mathrm l}
{1\mskip-4.5mu\mathrm l} {1\mskip-5mu\mathrm l}}}

 \newcommand{\kommentar}[1]{}

\renewcommand{\subsection}{\secdef \subsct\sbsect}
\newcommand{\subsct}[2][default]{\refstepcounter{subsection}
\vspace{0.15cm}
{\flushleft\bf \arabic{section}.\arabic{subsection}~\bf #1  }
\nopagebreak\nopagebreak}
\newcommand{\sbsect}[1]{\vspace{0.1cm}\noindent
{\bf #1}\vspace{0.1cm}}

{\nopagebreak {\hfill{$\diamond$}}\\ }

\newtheorem{theorem}{Theorem}[section]

\newtheorem{prop}[theorem]{Proposition}

\theoremstyle{definition}

\def\thebibliography#1{\section*{Bibliography}
  \list%
  {\arabic{enumi}.}
    {\settowidth\labelwidth{[#1]}\leftmargin\labelwidth
    \advance\leftmargin\labelsep
    \parsep0pt\itemsep0pt
    \usecounter{enumi}}
    \def\newblock{\hskip .11em plus .33em minus .07em}
    \sloppy                   
    \sfcode`\.=1000\relax}

\def\P{\prob}

\newcommand{\V}{\mathbb{V}}
\newenvironment{proofsect}[1]
{\vskip0.1cm\noindent{\bf #1.}}{\vspace{0.15cm}}

\begin{document}
\title[Local times of a transient random walk on $\Z^d$]{\Large Moments and distribution of the local times\\ of a transient random walk on $\Z^d$}
\author[Mathias Becker and Wolfgang K\"onig]{}
\maketitle
\thispagestyle{empty}

\centerline{\sc By Mathias Becker\footnote{Universit\"at Leipzig, Mathematisches Institut, Postfach 10 09 20, D-04009 Leipzig, Germany}$^{, }$\footnote{\tt mbecker@math.uni-leipzig.de} and Wolfgang K\"onig$^{1, }$\footnote{\tt koenig@math.uni-leipzig.de}}
\renewcommand{\thefootnote}{}
\footnote{\textit{AMS 2000 Subject Classification.} 60G50, 60J55, 60F15 }
\footnote{\textit{Keywords.} Random walk on $\Z^d$, local time, self-intersection number}
\renewcommand{\thefootnote}{1}

\vspace{0.2cm}
\centerline{\small \version}
\vspace{0.2cm}

\begin{quote}{\small }{\bf Abstract.} Consider an arbitrary transient random walk 
on $\Z^d$ with $d\in\N$. Pick $\alpha\in[0,\infty)$ and let $L_n(\alpha)$ be the spatial sum of the $\alpha$-th power of the $n$-step local times of the walk. Hence, $L_n(0)$ is the range, $L_n(1)=n+1$, and for integers $\alpha$, $L_n(\alpha)$ is the number of the $\alpha$-fold self-intersections of the walk. We prove a strong law of large numbers for $L_n(\alpha)$ as $n\to\infty$. Furthermore, we identify the asymptotic law of the local time in a random site uniformly distributed over the range. These results complement and contrast analogous results for recurrent walks in two dimensions recently derived by \v{C}ern\'y \cite{Ce07}. Although these assertions are certainly known to experts, we could find no proof in the literature in this generality.
\end{quote}


\section{Introduction and main results}

\noindent Let $(X_i)_{i\in\N}$ be a sequence of independent, identically distributed, $\Z^d$-valued random variables. Let $(S_n)_{n\in\N_0}$ be the corresponding random walk:
\begin{eqnarray}\label{gleichung01}
S_0:=0&\textrm{ and }&
S_n:=\sum_{i=1}^n X_i, \qquad n\in\N.
\end{eqnarray}
The main object of the present paper are the so-called {\it local times}
\begin{equation}
\ell(n,x):=\sum_{i=0}^n\1_{\{ S_i=x\} },\qquad n\in\N, x\in\Z^d,
\end{equation}
the number of visits to $x$ by time $n$. More specifically, we are interested in the large-$n$ asymptotics of the following functional of the local times:
\begin{eqnarray}
L_n(\alpha)&:=&\sum_{x\in\Z^d}\ell(n,x)^\alpha, \qquad \alpha\geq 0. 
\end{eqnarray}
This is a rather natural object in the study of random walks on $\Z^d$. Much attention has been focused on $L_n(0)=|\{S_0,\ldots,S_n\}|$, the {\it range} of the random walk, which is the number of distinct lattice points visited up to time $n$. The case $\alpha=1$ is trivial, as $L_n(1)=n+1$. Furthermore, $L_n(2)$ is the {\it self-intersection local time}, the number of self-intersections, which has been much studied from physical motives. More generally, for $\alpha$ an integer, $L_n(\alpha)$ is the number of $\alpha$-fold self-intersections of the random walk up to time $n$. The quantity $L_n(2)$ also arises as the variance of the {\it random walk in random scenery} \cite{KS79}.

In this paper, we first prove a strong law of large numbers for $L_n(\alpha)$ for any $\alpha\geq 0$, and we second show that $\ell(n,Y_n)$, where $Y_n$ is uniformly distributed on the set $\{S_0,\ldots,S_n\}$, has asymptotically a geometric distribution. We restrict to transient random walks in all dimensions. Our precise assumptions in this paper are the following.

We assume that the random walk is genuinely $d$-dimensional in the sense of \cite[D7.1]{Sp64}. That is, the set $R^+-R^+$ is $d$-dimensional, where $R^+=\bigcup_{n\in\N_0}\{x\in\Z^d\colon \P(S_n=x)>0\}$ is the support of the random walk. Furthermore, we assume that the {\it escape probability} of the random walk,
\begin{equation}
\gamma=\P(S_n\neq 0 \textrm{ for any }n\in\N),
\end{equation}
satisfies $0<\gamma<1$. (The condition that $\gamma<1$ rules out trivial cases). Finally, in the case $d\in\{1,2\}$, we assume that either the second moment of the steps is finite or that there is some $\eta>0$ such that $\sum_{k=n}^\infty \P(S_k=0)\leq O(n^{-\eta})$ as $n\to\infty$.

We have not been able to decide if any of our two conditions in the case $d\in\{1,2\}$ implies the other or not. Now we formulate the strong law of large numbers. 
\begin{theorem}\label{theorem1}
For all $\alpha\in[0,\infty)$ it holds $\P$-almost surely:
\begin{eqnarray}\label{gleichung1}
\lim_{n \to \infty} \frac{L_n(\alpha )}{n} 
& = & \sum_{j\in\N}j^\alpha \gamma^2(1-\gamma)^{j-1}.
\end{eqnarray}
\end{theorem}

Note that the assertion for $\alpha=0$ is already well-known; a proof is contained, e.g., in \cite[T4.1]{Sp64}. 

\begin{theorem}\label{theorem2}
Given the steps $X_1,\ldots,X_n$, let $Y_n$ be a uniformly distributed random variable on the set of visited points, $\{ S_0,\ldots ,S_n\}$. Then, $\P$-almost surely,
\begin{eqnarray}
\lim_{n\to\infty}\mathbb{P}\left[\left. \ell(n,Y_n)= u\right| X_1,\ldots ,X_n\right]&=& \gamma(1-\gamma)^{u-1}, \qquad u\in\N.
\end{eqnarray}
\end{theorem}
We prove Theorems \ref{theorem1} and \ref{theorem2} in Section \ref{sect_thmproof}. In Section \ref{sect_expval} we analyse the expected value of $L_n(\alpha)$ and in Section \ref{sect_var} its variance.

Theorem \ref{theorem2} can be understood as follows. For large $n$, the event $\{\ell(n,Y_n)=u\}$ is realised by returning $u-1$ times to $Y_n$, which has probability approximately $1-\gamma$ each, and not returning afterwards, which has probability $\approx \gamma$. This makes plausible that $\ell(n,Y_n)$ converges in distribution towards a geometric random variable $Z$ with parameter $\gamma$. Observe that the right-hand side of \eqref{gleichung1} is equal to $\gamma\E(Z^\alpha)$. Hence also Theorem \ref{theorem1} can easily be understood, taking into consideration the well-known fact that $|\{S_0,\ldots,S_n\}|\sim \gamma n$ (see \cite[T4.1]{Sp64}), which means that the sum on $x$ in the definition of $L_n(\alpha)$ has effectively only $\approx \gamma n$ summands.

It is remarkable that Theorems \ref{theorem1} and \ref{theorem2} hold for an arbitrary, genuinely $d$-dimensional random walk in $d\geq 3$, without any integrability, centering or periodicity condition. This is in sharp contrast to the recurrent, two-dimensional case studied in \cite{Ce07}, where it was assumed that the steps are centred and have finite second moments. Also observe that our results contain the low-dimensional case with non-zero drift, which leads to the same asymptotics as in the higher-dimensional case.

Analogous results for the two-dimensional case have recently been derived in \cite{Ce07}. 
To the best of our knowledge Theorem \ref{theorem2} has not yet appeared in the literature, although related assertions already appeared in \cite{ET60}. Various special cases of Theorem~\ref{theorem1} are spread over the literature. \cite{DE51} considered the asymptotics of $\E(L_n(0))$ for simple random walk. An extension to a somewhat more general random walk can be found in \cite{We94}. Some estimates of the variance of $L_n(2)$ for simple random walk have been derived in \cite{BS95}. The novelty and main value of the present paper lies in its generality in three aspects: all $\alpha\geq 0$, very general transient random walks, and the strong law of large numbers.

\section{The expected value $\E(L_n(\alpha))$} \label{sect_expval}

\noindent In this section we prove the asymptotics in (\ref{gleichung1}) for the expected value of $L_n(\alpha)$ in place of the variable $L_n(\alpha)$ itself. We use the approach of \cite{ET60} (which considers simple random walk only) and use the opportunity to correct an error in \cite[Theorem 12]{ET60}.

\begin{prop}\label{satz35}
For any $\alpha \in[0,\infty)$,
\begin{eqnarray}
\lim_{n\to\infty}\frac{\E (L_n(\alpha))}{n}
&=&\sum_{j\in\N}j^\alpha\gamma^2 (1-\gamma)^{j-1}.
\end{eqnarray}
\end{prop}
\begin{proofsect}{Proof} 
Introduce the probability that the point visited by random walk in the $n$-th step has not been visited before:
$$
\gamma(n)=\P(S_0\neq S_n,\ldots,S_{n-1}\neq S_n),
$$
in particular $\gamma(0)=1$. Note that $1-\gamma(1)=\P(S_1=0)$ may be zero.
An easy computation yields that $\gamma(n)$ is also equal to the probability that the random walk does not return to the origin within the first $n$ steps:
$$
\gamma(n)=\P(S_1\neq S_0,\ldots,S_{n}\neq S_0).
$$
Let $\tau=\inf\{n\in\N\colon S_0=S_n\}$ denote the return time to the origin, then we have $\P(\tau=n)=\gamma(n-1)-\gamma(n)$ for any $n\in\N$. 
By the monotone convergence theorem, $\gamma=\lim_{n\to\infty}\gamma(n)$.

We introduce now the number of points that have been visited exactly $j$ times up to time $n$:
$$
Q_j(n)=|\{x\in\Z^d\colon \ell(n,x)=j\}|. 
$$
Its expectation can be calculated as follows.
\begin{equation}\label{gleichung33}
\begin{aligned}
\E\left(Q_j(n)\right)&=\sum_{x\in\Z^d}\P(\ell(n,x)=j)\\
&=\sum_{\substack{{x\in\Z^d}\\ {0\leq k_1<k_2<\ldots <k_j\leq n}}}\P\big(S_{k_1}=\dots=S_{k_j}=x, S_k\not=x\mbox{ for }k\in\{0,\dots,n\}\setminus\{k_1,\dots,k_j\}\big)\\
&=\sum_{0\leq k_1<k_2<\ldots <k_j\leq n}\gamma(k_1)\left[\prod_{i=1}^{j-1}\P(\tau=k_{i+1}-k_{i})\right]\gamma(n-k_j),
\end{aligned}
\end{equation}
where we used the Markov property at the times $k_1,\dots,k_j$ and the definition of $\tau$. 
We use (\ref{gleichung33}) to identify the generating function $q_j(s)=\sum_{n\in\N_0}\E\left(Q_j(n)\right)s^n$ of $(Q_j(n))_{n\in\N}$ as
\begin{equation}\label{gleichung51}
q_j(s)=\left(\sum_{n=0}^\infty s^{n}\gamma(n)\right)^2\left( \sum_{n=1}^\infty s^n \P(\tau=n)\right)^{j-1}.
\end{equation}
Since $\lim_{n\to\infty}\gamma(n)=\gamma$ and by the monotone convergence theorem, its behaviour as $s\nearrow 1$ is identified as follows
\begin{equation}\label{qjasy}
q_j(s)\sim\left(\frac{\gamma}{1-s}\right)^2\left( \sum_{n=1}^\infty s^n \P(\tau=n)\right)^{j-1} \sim(1-s)^{-2}\gamma^2(1-\gamma)^{j-1},
\end{equation}
where we note that $\sum_{n=1}^\infty \P(\tau=n)=\P(\tau<\infty)=1-\gamma$. 
By the Tauberian theorem (see \cite[Theorem XIII.5]{Fe71}) we obtain
\begin{equation}\label{gleichung31}
\lim_{n\to\infty}\E\left(\frac{Q_j(n)}{n}\right)= \gamma^2\left( 1-\gamma\right)^{j-1}.
\end{equation}
(At this point we see that in \cite[Theorem 12]{ET60} a factor $\gamma$ is missing. Let us remark that \eqref{gleichung31} was also derived in \cite{P74} using entirely different methods.)

Now we return to the expected value of $L_n(\alpha)$. It is easy to see that
$$
\E (L_n(\alpha))=\sum_{j\in\N}j^\alpha\E\left(Q_j(n)\right).
$$
We will now give an upper and a lower asymptotic bound for the generating function of $(\E(L_n(\alpha)))_{n\in\N}$ which will turn out to be identical. With the help of (\ref{gleichung51}), similarly to \eqref{qjasy}, we deduce, for $s\in [0,1)$,
$$
\begin{aligned}
\sum_{n\in\N}s^n\E(L_n(\alpha))=\sum_{j\in\N}j^\alpha q_j(s)
&\leq \sum_{j\in\N}j^\alpha\left(\sum_{n=0}^\infty s^{n}\gamma(n)\right)^2\left(1-\gamma\right)^{j-1}\\
&\sim\left(\frac{\gamma}{1-s}\right)^2\sum_{j\in\N}j^\alpha\left(1-\gamma\right)^{j-1},\qquad\mbox{as }s\nearrow 1.
\end{aligned}
$$
For the lower bound we use Fatou's lemma and obtain:
$$
\begin{aligned}
\liminf_{s\uparrow 1}(1-s)^2\sum_{n\in\N}s^n\E(L_n(\alpha))
&=\liminf_{s\uparrow 1}\sum_{j\in\N}j^\alpha(1-s)^2 q_j(s)\\
&\geq\sum_{j\in\N}j^\alpha\liminf_{s\uparrow 1}(1-s)^2 q_j(s)
=\sum_{j\in\N}j^\alpha\gamma^2(1-\gamma)^{j-1}.
\end{aligned}
$$
Now apply once more the Tauberian theorem to complete the proof of (\ref{satz35}).
\end{proofsect}
\qed

\section{The variance $\V(L_n(\alpha))$}\label{sect_var}

\noindent In this section we prove an upper bound on the variance of $L_n(\alpha)$ that is sufficient for the application of the second-moment method in Section \ref{sect_thmproof}. For this, it suffices to show that $\V(L_n(\alpha))=o(\E(L_n(\alpha))^{2-\kappa})=o(n^{2-\kappa})$ for some $\kappa>0$ (see Proposition~\ref{satz35}). We now give a bound that seems optimal in $d\geq 3$. Our method is a combination of ideas from \cite{Ce07} (which are based on ideas from \cite{Bo89}) and from the proof of \cite[Prop.~3.1]{BS95}, where the normalized self-intersection number of simple random walk is estimated. We are also able to apply a result from \cite{JP71}.

\begin{prop} \label{theorem32}
Fix $\alpha\in\N$. 
\begin{enumerate} 
\item[(i)] In the case $d\in\{1,2\}$, if the walker's steps have a finite second moment, then there exists a constant $C>0$ such that, for any $n\in\N$,
\begin{equation}\label{variancebounda}
\V (L_n(\alpha))\leq C\times\begin{cases} n^{3/2}\log n & \textrm{in }d=1,
 \\ n\log^2 n & \textrm{in }d=2.
 \end{cases}
\end{equation}
On the other hand, if $\sum_{k=n}^\infty \P(S_k=0)=O(n^{-\eta})$ for some $\eta>0$, then
there exists a constant $C>0$ such that,
\begin{equation}\label{varianceboundb}
\V (L_n(\alpha))\leq C n^{2-\eta},\qquad n\in\N.
\end{equation}

\item[(ii)] In the case $d\geq 3$, there exists a constant $C>0$ such that, for any $n\in\N$,
\begin{equation}\label{variancebound}
\V (L_n(\alpha))\leq C\times\begin{cases} n^{3/2} & \textrm{in } d=3,
 \\ n\log n & \textrm{in }d=4, 
 \\ n & \textrm{in }d\geq 5.
 \end{cases}
\end{equation}
\end{enumerate}
\end{prop}

\begin{proofsect}{Proof} Since $\alpha$ is an integer we may rewrite the variance as follows.
$$
\begin{aligned}
\V (L_n(\alpha )) 
&= \V \left[ \sum_{x\in\mathbb{Z}^d} \sum_{k_1,\ldots ,k_\alpha =0}^n \1_{\{ S_{k_1}=\ldots =S_{k_\alpha}=x\}}\right]\\
&= \sum_{\heap{k_1,\ldots ,k_\alpha =0} {l_1,\ldots ,l_\alpha =0}}^n \big[\P (S_{k_1}=\ldots =S_{k_\alpha}, S_{l_1}=\ldots =S_{l_\alpha} )-\P (S_{k_1}=\ldots =S_{k_\alpha})\P (S_{l_1}=\ldots =S_{l_\alpha} )\big]\\
&= \sum_{\beta, \gamma =1} ^\alpha C(\alpha ,\beta ,\gamma ) a_{\beta ,\gamma}(n),
\end{aligned}
$$
where 
\begin{equation}\label{andef}
a_{\beta ,\gamma}(n)=\sum_{\heap{0\leq k_1<\ldots <k_\beta\leq n}{ 0\leq l_1<\ldots <l_\gamma\leq n}}
\big[\P (S_{k_1}=\ldots =S_{k_\beta}, S_{l_1}=\ldots =S_{l_\gamma} )-\P (S_{k_1}=\ldots =S_{k_\beta})\P (S_{l_1}=\ldots =S_{l_\gamma} ) \big]
\end{equation}
and
$$
\begin{aligned}
C(\alpha,\beta,\gamma)
&=\big|\big\{(\widetilde k_1,\ldots,\widetilde k_\alpha)\in\{0,\ldots,n\}^\alpha\colon \{k_1,\ldots,k_\beta\}=\{\widetilde k_1,\ldots,\widetilde k_\alpha\}\big\}\big|\\
&\quad\times\big|\big\{(\widetilde l_1,\ldots,\widetilde l_\alpha)\in\{0,\ldots,n\}^\alpha\colon \{l_1,\ldots,l_\gamma\}=\{\widetilde l_1,\ldots,\widetilde l_\alpha\}\big\}\big|.
\end{aligned}
$$
In words: $C(\alpha ,\beta ,\gamma )$ is the number of pairs of unordered tuples $(\widetilde k_1,\ldots ,\widetilde k_\alpha)$ and $(\widetilde l_1,\ldots ,\widetilde l_\alpha)$ (with possible repetitions) that give the same sequence $0\leq k_1<\ldots <k_\beta\leq n$ respectively $0\leq l_1<\ldots <l_\gamma\leq n$. Note that $C(\alpha ,\beta ,\gamma )$ only depends on $\alpha,\beta$ and $\gamma$. Hence, we only have to show that $a_{\beta ,\gamma}(n)$ satisfies the bound in \eqref{variancebound} for any $\beta$ and $\gamma$.

Let $\beta,\gamma\in\{1,\dots,\alpha\}$ be fixed. In order to prepare for the application of the Markov property, we encode the tuples of the numbers $0\leq k_1<\ldots <k_\beta\leq n$ and $ 0\leq l_1<\ldots <l_\gamma\leq n$ in terms of a new set of variables $(j_i,\kappa_i ) \in \{ 0,\ldots ,n\} \times  \{ 0,1\}$ for $i\in\{1,\ldots,\beta+\gamma\}$ satisfying
$$
 \{ j_i\colon \kappa_i =0 \} =\{ k_1,\ldots ,k_\beta\} \qquad\mbox{and}\qquad
 \{ j_i\colon \kappa_i =1 \} =\{ l_1,\ldots ,l_\gamma\},
$$
and such that $(j_i,\kappa_i)_{i=1}^{\beta+\gamma}$ is alphabetically ordered, i.e.,  $j_i \leq  j_{i+1}$ for any $i$, and if $j_i=j_{i+1}$ then $\kappa_i<\kappa_{i+1}$. Then we introduce $m_0=j_1$ and $m_{\beta +\gamma} = n-j_{\beta +\gamma}$ and
$$
m_i =j_{i+1}-j_{i}\qquad\mbox{and}\qquad \eps_i \ = \kappa_{i+1}-\kappa_i\qquad\mbox{for}\quad 
i=1,\dots,\beta +\gamma-1.
$$
In this way, we have mapped tuples of numbers $0\leq k_1<\ldots <k_\beta\leq n$ and $ 0\leq l_1<\ldots <l_\gamma\leq n$ one-to-one onto sequences $(m_i,\eps_i)_{i=0}^{\beta+\gamma}$ in $\N_0\times\{-1,0,1\}$ satisfying $\sum_{i=0}^{\beta+\gamma}m_i=n$. In words, $m_i$ is the difference between the $(i+1)$-th largest and the $i$-th largest of the numbers $k_1,\dots,k_\beta,l_1,\dots,l_\gamma$, and $\eps_i$ is 1 if the $m$-sequence switches from an $k$-value to an $l$-value, and $\eps_i$ is $-1$ if it switches in the reversed way. Since $\beta\geq 1$ and $\gamma\geq 1$, one easily sees that $\#\{i\colon\eps_i\neq 0\}\geq 1$. 

Let us look at the case that $\#\{i\colon\eps_i\neq 0\}$ is equal to one, i.e., either the $k_j$'s all all smaller than the $l_j$'s or the other way around. By use of the Markov property, we easily see that in this case all the summands on the right-hand side of the definition \eqref{andef} of  $a_{\beta ,\gamma}(n)$ vanish. Thus we may restrict to the case $\#\{i\colon\eps_i\neq 0\}\geq 2$. In this case, it is clear that $\eps_u=1=-\eps_v$, for some $u,v\in\{1,\dots,\beta+\gamma-1\}$.
For our purposes, it will turn out to suffice to estimate the negative terms on the right-hand side of \eqref{andef} against zero. Hence, we obtain
$$
\begin{aligned}
a_{\beta ,\gamma}(n)
&\leq \sum_{\substack{{0\leq k_1<\ldots <k_\beta\leq n}\\ {0\leq l_1<\ldots <l_\gamma\leq n}}} \1_{\{\exists u,v\colon \eps_u=1=-\eps_v\}}
\P (S_{k_1}=\ldots =S_{k_\beta}, S_{l_1}=\ldots =S_{l_\gamma} )\\
&\leq    \sum_{\substack{{m_0, m_1, \ldots, m_{\beta+\gamma}\in\N_0}\\ {m_0+m_1+\ldots+m_{\beta+\gamma}=n}\\ {\eps\in\{-1,0,1\}^{\beta+\gamma+1}}}} \sum_{u,v=1}^{\beta+\gamma-1}\1_{\{\eps_u=1=-\eps_v\}}\sum_{x,y\in\Z^d}
\P (S_{m_0}=x)\prod_{i=1}^{\beta-\gamma-1}\P(S_{m_i}=\eps_i y),
\end{aligned}
$$
where $x$ plays the role of $S_{k_1\wedge l_1}$ and $y$ the role of $S_{l_1}-S_{k_1}$ respectively $S_{k_1}-S_{l_1}$. 

In the following, we use $C$ to denote a generic positive constant that depends on $\beta$, $\gamma$, $d$ and the step distribution only and may change its value from appearance to appearance. The summation over $m_{\beta+\gamma}$ is redundant, as $m_{\beta+\gamma}$ can be computed from the other $m_i$ and $n$. Furthermore, we execute the sum of the terms $\P (S_{m_0}=x)$ over $x\in\Z^d$, and consequently the summation over $m_0$ delivers an additional factor $n$ as an upper bound. Now we also execute the sum over all $m_i$ with $i\not= u,v$ and use that $\sup_{y\in\Z^d}\sum_{m_i=0}^\infty \P(S_{m_i}=\eps_i y)\leq C$, by transience. Hence, we obtain
\begin{equation}\label{equa2}
\begin{aligned}
a_{\beta ,\gamma}(n)
&\leq Cn\sum_{\substack{{u,v=1}\\{u\not= v}}}^{\beta+\gamma-1}\sum_{y\in\Z^d}\sum_{m_u,m_v=0}^n \P(S_{m_u}=y)\P(S_{m_v}=-y)\\
&\leq Cn \sum_{y\in\Z^d}\sum_{m,\widetilde m=0}^n \P(S_{m}=y)\P(S_{\widetilde m}=-y).
\end{aligned}
\end{equation}

Now, in dimensions $d\geq 3$, our assertion in \eqref{variancebound} directly follows from \cite[Lemma 3]{JP71}, which implies that, for every $n\in\N$,
\begin{equation}\label{JP}
\sum_{y\in\Z^d}G_n(0,y)G_n(0,-y)\leq C\times\begin{cases}n^{1/2}&\mbox{for }d=3,\\
\log n&\mbox{for }d=4,\\
1&\mbox{for }d\geq 5,
\end{cases}
\end{equation}
where $G_n(0,y)=\sum_{m=1}^n\P(S_{m}=y)$ denotes the $n$-step Green's function. The main tool in the proof of \eqref{JP} is the estimate
\begin{equation}\label{Spitzer}
\sup_{x\in\Z^d}\P(S_m=x)\leq Cm^{-\frac{d}{2}},\qquad m\in\N,
\end{equation}
which is proved in any dimension $d$ as \cite[Lemma 1]{JP71}, based on \cite[P7.6]{Sp64}. Hence, the proof of (ii) is finished.

Now we consider the case $d\in\{1,2\}$. First we work under the assumption of the existence of the second moment of the steps. Consider their expected value, $v=\E(X_1)\in\R^d\setminus\{0\}$. Consider the closed half-plane, respectively half-axis, $A=\{x\in\R^d\colon x\cdot v\leq 0\}$, where $x\cdot v$ is the standard inner product. Then $-A^{\rm c}\subset A$. Splitting the sum on the right hand side of \eqref{equa2} into $y\in A$ and $y\in A^{\rm c}$ and using the symmetry of the summands in $y$, we obtain:  
\begin{equation}\label{eqn3}
\begin{aligned}
a_{\beta ,\gamma}(n)
&\leq  Cn\sum_{y\in A}\sum_{m,{\widetilde m}=1 }^n\P(S_{m}=y)\P(S_{{\widetilde m}}=-y)\\
&\leq  Cn\left(\sum_{m=1}^n\sum_{y\in A}\P(S_{m}=y)\right)\sum_{{\widetilde m}=1}^n\sup_{y\in\Z^d}\P(S_{{\widetilde m}}=-y)\leq Cn\sum_{m=1}^n\P(S_m\in A)\sum_{{\widetilde m}=1}^n\widetilde m^{-\frac d2}, 
\end{aligned}
\end{equation} 
where we also used \eqref{Spitzer}. We have, for any $m\in\N$, using Chebyshev's inequality,
\begin{equation}
\P(S_{m}\in A)\leq \P(|S_{m}-vm|\geq \|v\| m)=\P(|S_{m}-\E(S_m)|\geq\|v\| m)\leq \frac{\V(S_{m})}{\|v\|^2 m^2}\leq \frac{C}{m}.
\end{equation}
Using this in \eqref{eqn3} yields the result in \eqref{variancebounda}.

Finally, we work in $d\in\{1,2\}$ under the assumption that $\sum_{k=n}^\infty \P(S_k=0)=O(n^{-\eta})$. We go back to \eqref{equa2} and use the Markov property to rewrite
$$
a_{\beta ,\gamma}(n)\leq Cn \sum_{m,\widetilde m=0}^n \P(S_{m+\widetilde m}=0)
\leq Cn\sum_{j=0}^{2n}\sum_{k=j}^\infty \P(S_k=0)\leq Cn \sum_{j=0}^{2n}j^{-\eta}\leq Cn^{2-\eta}.
$$
This means that we have derived \eqref{varianceboundb}. The proof of the proposition is complete.
\end{proofsect}
\qed

\section{Proof of the theorems}\label{sect_thmproof}

\noindent In this section we prove our two main results, Theorems~\ref{theorem1} and \ref{theorem2}. Our strategy is the same as in \cite{Ce07}; the proof will be done in three steps. First we prove the assertion in Theorem \ref{theorem1} for all integers $\alpha\geq 0$. Second we provide the proof for Theorem \ref{theorem2}, and the last step is the proof of the general assertion in Theorem \ref{theorem1} for $\alpha\in[0,\infty)$.

\begin{proofsect}{Proof of Theorem \ref{theorem1} for $\boldsymbol{\alpha\in\N_0}$} For $\alpha=0$ the assertion is proved in \cite[T4.1]{Sp64}. For $\alpha=1$ the assertion is obvious, so let $\alpha\geq2$. From Proposition~\ref{theorem32} we have, for some $\kappa>0$, that $\V(L_n(\alpha))\leq C n^{2-\kappa}$ for all $n\in\N$, with a suitable constant $C>0$. Consider the subsequence $n_k:=\lfloor k^{2/\kappa}+1\rfloor$. With the help of Chebyshev's inequality, Propositions~\ref{satz35} and \ref{theorem32}, we get, for a suitable constant $C$ and any $\eps>0$:
\begin{equation}
\P\big(|L_{n_k}(\alpha)-\E L_{n_k}(\alpha)|\geq\eps\E (L_{n_k}(\alpha))\big)
\leq \frac{\V(L_{n_k}(\alpha))}{\eps^2\left(\E (L_{n_k}(\alpha))\right)^2}\leq C\eps^{-2}n_k^{-\kappa}\leq \frac{C}{\eps^2k^2}.
\end{equation}
Since this is summable over $k\in\N$, the Borel-Cantelli lemma yields that, $\P$-almost surely,
\begin{equation}
\lim_{k\to\infty}\frac{L_{n_k}(\alpha)}{\E(L_{n_k}(\alpha))} =1. 
\end{equation}
It remains to fill the gaps between the $n_k$. With the help of the monotonicity of $n\mapsto L_n(\alpha)$ and the fact that $\lim_{k\to\infty}\E(L_{n_{k+1}}(\alpha))/\E(L_{n_k}(\alpha))=\lim_{k\to\infty}n_{k+1}/n_k=1$ (see Proposition~\ref{satz35}), we can easily deduce that, $\P$-almost surely,
$$
\lim_{n\to\infty}\frac{L_n(\alpha)}{\E(L_n(\alpha))}=1.
$$
Now Proposition~\ref{satz35} finishes the proof of Theorem \ref{theorem1} for all $\alpha\in\N$.
\end{proofsect}
\qed

\begin{proofsect}{Proof of Theorem \ref{theorem2}} Let $Y_n$ be uniformly distributed on $\{S_0,\ldots,S_n\}$, and let $Z_n = \ell(n,Y_n)$. Denote $R(n)=|\{S_0,\ldots,S_n\}|=L_n(0)$ and recall that $\lim_{n\to\infty}R(n)/n=\gamma$, according to \cite[T4.1]{Sp64}. Hence, for $\alpha\in\N_0$, the conditional $\alpha$-th moments of $Z_n$ given $X_1,\ldots,X_n$ can be asymptotically identified as
\begin{eqnarray}\label{gleichung37}
\E(Z_n^\alpha|X_1,\ldots,X_n) 
= \sum_{x\in \Z^d}\frac{1}{R(n)}\ell(n,x)^\alpha
= \frac{\gamma n}{R(n)}\frac{L_n(\alpha)}{\gamma n} 
\sim\sum_{j\in\N}j^\alpha\gamma (1-\gamma)^{j-1}.
\end{eqnarray}
In other words, these moments converge to the moments of a geometrically distributed random variable with parameter $\gamma$. With the help of the theorem of Fr\'echet-Shohat (see \cite[V.1]{Sc98}) we get the weak convergence of the conditional distribution of $Z_n$ given $X_1,\ldots,X_n$ towards a geometric distribution with the parameter $\gamma$. This finishes the proof of Theorem \ref{theorem2}.
\end{proofsect}
\qed

\begin{proofsect}{Proof of Theorem \ref{theorem1} in the general case} Recall that it suffices to consider $\alpha\in[0,\infty)\setminus\N_0$. The continuous mapping theorem and Theorem \ref{theorem2} imply the weak convergence of the conditional distribution of $Z_n^\alpha$ given $X_1,\ldots,X_n$ towards $Z^\alpha$ if $Z$ is geometrically distributed with parameter $\gamma$. Since that sequence is uniformly integrable, we also have  
\begin{equation}\label{gleichung39}
\lim_{n\to\infty}\E\left(Z_n^\alpha|X_1,\ldots,X_n\right) =\E(Z^\alpha)
=\sum_{j\in\N}j^\alpha\gamma (1-\gamma)^{j-1}. 
\end{equation}
As in (\ref{gleichung37}), we see that 
\begin{equation}
\frac{L_n(\alpha)}{n}
=\frac{R(n)}{n}\E\left(Z_n^\alpha|X_1,\ldots,X_n\right).
\end{equation}
Now use \eqref{gleichung39} and the fact that $\lim_{n\to\infty}R(n)/n=\gamma$ (recall \cite[T4.1]{Sp64}) to finish the proof of Theorem \ref{theorem1}.
\end{proofsect}
\qed

{\bf Acknowledgement.} We thank two anonymous reviewers whose suggestions lead to improvements of our results.



\end{document}